\documentclass[11pt,reqno]{amsart}

\usepackage{amsmath, amsfonts, amssymb, amscd, enumerate}
\theoremstyle{plain}
\newtheorem{theo}{Theorem}[section]
\newtheorem{lem}[theo]{Lemma}
\newtheorem{prop}[theo]{Proposition}

\theoremstyle{definition}
\newtheorem{rem}[theo]{Remark}

\newtheorem{definition}[theo]{Definition}
\newenvironment{pf}{\noindent{\it Proof. }}{$\square$\par\medskip}

\theoremstyle{plain}

\theoremstyle{definition}

\renewcommand{\=}{\overset{\operatorname{def}}{=}}

\newcommand{\beq}{\begin{equation}}
\newcommand{\eeq}{\end{equation}}
\renewcommand{\a}{\alpha}
\renewcommand{\b}{\beta}

\renewcommand{\d}{\delta}

\newcommand{\g}{\gamma}

\renewcommand{\l}{\lambda}

\newcommand{\s}{\sigma}

\newcommand{\z}{\zeta}
\newcommand{\D}{\Delta}

\newcommand{\G}{\Gamma}

\newcommand{\bC}{\mathbb{C}}
\newcommand{\bR}{\mathbb{R}}

\newcommand{\bJ}{\mathbb{J}}

\newcommand{\bN}{\mathbb{N}}

\newcommand{\gR}{\mathfrak{R}}

\newcommand\GL{\mathrm{GL}}

\newcommand{\cC}{\mathcal{C}}
\newcommand{\cD}{\mathcal{D}}

\newcommand{\cF}{\mathcal{F}}
\newcommand{\cG}{\mathcal{G}}

\newcommand{\cL}{\mathcal{L}}

\newcommand{\cN}{\mathcal{N}}

\newcommand{\cR}{\mathcal{R}}

\newcommand{\cU}{\mathcal{U}}
\newcommand{\cV}{\mathcal{V}}
\newcommand{\cW}{\mathcal{W}}

\newcommand{\Jst}{J_{\operatorname{st}}}

\renewcommand{\square}{\kern1pt\vbox
{\hrule height 0.6pt\hbox{\vrule width 0.6pt\hskip 3pt
\vbox{\vskip 6pt}\hskip 3pt\vrule width 0.6pt}\hrule height0.6pt}\kern1pt}

\DeclareMathOperator\Aut{Aut\;}

\renewcommand\Re{\operatorname{Re}}
\renewcommand\Im{\operatorname{Im}}

\newcommand{\wt}{\widetilde}
\newcommand{\wh}{\widehat}

\newcommand{\n}{\nabla}

\newcommand{\dist}{\operatorname{dist}}

\newcommand{\be}{\begin{equation}}
\newcommand{\ee}{\end{equation}}

\def\<#1,#2>{\langle\,#1,\,#2\,\rangle}
\newcommand{\arr}{\begin{array}{rlll}}
\newcommand{\ea}{\end{array}}
\newcommand{\bea}{\begin{eqnarray}}
\newcommand{\eea}{\end{eqnarray}}
\newcommand{\bean}{\begin{eqnarray*}}
\newcommand{\eean}{\end{eqnarray*}}

\def\sideremark#1{\ifvmode\leavevmode\fi\vadjust{%            The remark
\vbox to0pt{\hbox to 0pt{\hskip\hsize\hskip1em%               will appear only
\vbox{\hsize3cm\tiny\raggedright\pretolerance10000%          on the side
\noindent #1\hfill}\hss}\vbox to8pt{\vfil}\vss}}}%           in 3cm
\newcounter{ssig}
\newcounter{ttig}
\newcommand{\upG}[1]{{#1}^{G}}

\title[Foliations  by stationary disks of  almost complex domains]
{Foliations by stationary disks\\ of  almost complex domains}

\author{G. Patrizio and A. Spiro}

\begin{document}

\begin{abstract}   We study the problem of existence of stationary disks for domains
in almost complex manifolds. As a consequence of our results, we prove that any almost complex domains which is a small deformations of a strictly linearly convex domain $D \subset \bC^n$ with standard complex structure admits a singular foliation by stationary disks passing through any given internal point. Similar results are given for foliation by stationary disks through a given boundary point.\par
\vskip 0.3 truecm
\noindent
\font\smallsmc = cmcsc9
 {\smallsmc R\'esum\'e}.
Nous \'etudions   le probl\`eme de l'existence des disques stationnaires pour des domaines dans une vari\'et\'e presque complexe. Comme cons\'equence de nos r\'esultats, nous montrons que tous les domaines presque complexes obtenu comme une petite 
d\'eformation  d'un domaine strictement lin\'eairement convexe 
$D \subset \bC^n$, avec la structure complexe standard, admet une foliation dans  disques stationnaires passant par un point interne 
donn\'e de $D$. Des r\'esultats similaires sont obtenus pour  foliations dans disques  stationnaires  dont les bords passent pour un point donn\'e dans le bord du domaine.
\end{abstract}

%\begin{abstract} We prove that for any point $x_o$ of a smoothly bounded, strictly linearly convex domain $D \subset \bC^n$, endowed with an almost complex structure   $J$,  sufficiently close to the standard complex structure$\Jst$, there exists a canonical singular foliation, formed by stationary disks passing through $x_o$. The result extends a theorem by Coupet, Gaussier and Sukhov  on the unit ball $B^n \subset \bC^n$, endowed with a small deformation $J$ of $\Jst$,  and establishes  a canonical diffeomorphism between  $(B^n, \Jst)$ and $(D, J)$ (the so-called ``generalized Riemann map'') for a large class of almost complex domains.  For the same  class of almost complex domains,  we also   show the existence of  regular foliations of  conical sectors made by  stationary disks passing through a fixed boundary point.  
% \end{abstract}


\thanks{{\it Acknowledgments}. This research was partially supported by the Project MIUR ÒGeometric Properties of Real and Complex ManifoldsÓ,  Project MIUR ÒDifferential Geometry and Global AnalysisÓ and by GNSAGA of INdAM}

\address{
%\vskip - 0.6cm 
Dipartimento di Matematica ``U. Dini'', 
Universit\`a di Firenze, 
%Viale Morgagni 67/a, 
Firenze, 
ITALY}
\email{patrizio@math.unifi.it}
\address{
%\vskip - 0.6cm 
Dipartimento di Matematica e Informatica, Universit\`a di Camerino, 
%Via Madonna delle Carceri, 
 Camerino, 
ITALY}
\email{andrea.spiro@unicam.it}
\maketitle

\null \vspace*{-.25in}

\section*{Introduction}
\setcounter{section}{1}
\setcounter{equation}{0}
Let $(M, J)$ be an almost complex manifold and $D \subset M$  a strongly  pseudoconvex domain with smooth boundary.  
Given a point $x_o \in D$, let us call  {\it foliation by stationary disks of $(D, x_o)$\/}  any  collection of  stationary disks centered  at $x_o$ and  smoothly parameterized by    the points of a unit sphere $S = \{ \ v \in T_{x_o} D\ :\ \|v\| = 1\}$ for  some Euclidean norm  $\|\cdot\|$ on $T_{x_o} D$. 
By ``stationary disk'' we mean any $J$-holomorphic embedding $f: \D \to D$ of the unit disk $\D \subset \bC$ that satisfies the definition of Coupet, Gaussier and Sukhov in \cite{CGS}, which  naturally generalizes  the usual notion of Lempert's stationary disks for bounded domains in $\bC^n$. \par
\smallskip
In case $(M,  J) = (\bC^n, \Jst)$,  natural examples of foliations by stationary disks are given  by the straight disks through  the origin of the pseudoconvex, smoothly bounded complete circular domains $D$ in $\bC^n$. Other interesting examples are provided by   the celebrated results by Lempert on  Kobayashi extremal disks  in  strictly linearly convex domains (\cite{Le, Le1, Le2}). In fact,  an immediate consequence of those results is  that for any  smoothly bounded, strictly linearly convex  domain $D \subset \bC^n$  and  any $x_o \in D$, the Kobayashi extremal disks of $D$ through $x_o$ give a foliation  by stationary disks of $(D, x_o)$. The existence of a foliation by stationary disks  is also one of the  main   properties of the smoothly bounded {\it domains of circular type}, a class of domains  in $\bC^n$ with an exhaustion of a special kind,  which naturally include all complete strictly pseudoconvex bounded circular domains, all bounded strictly linearly convex domains  and, more generally,  all strictly pseudoconvex domains with  (singular) foliations given  Kobayashi extremal disks satisfying some special regularity conditions (\cite{Pt1, Pt3}).\par
\medskip
In all these cases, the  foliation by stationary disks $\cF^{(x_o)}$  can be used to construct  a so-called 
 {\it generalized Riemann map\/}, i.e. 
a homeomorphism   $\varphi: \overline{B^n} \to \overline{D}$, 
which is smooth on $B^n \setminus \{0\}$ and maps the   straight complex lines in $B^n$ through $0$ 
 into corresponding disks of $\cF^{(x_o)}$.  This  generalized Riemann map have been often 
 used  in at least two important  research areas: a) generalizations of Fefferman's theorem on boundary regularity  of biholomorphisms between pseudoconvex domains; b) Green 
 functions with logarithmic pole for  Monge-Amp\`ere equations and plurisubharmonic exhaustions of pseudoconvex domains (see e.g. \cite{Le, Le1, Be, Tu, De}).\par  
\smallskip
At the best of our knowledge,  the first use of   foliations by stationary disks in the contest of  almost complex manifolds can 
 be found in \cite{CGS}.  There,  the authors generalize Lempert's notion of  stationary disks  in the almost complex  setting and show 
  the existence of a  foliation by  stationary disks of  the unit ball $B^n \subset \bC^n$, endowed with an almost complex structure $J$ which is a sufficiently small deformation of  standard complex structure $\Jst$. The corresponding 
 generalized Riemann map has been used  to prove  $\cC^\infty$-regularity of biholomorphisms between two almost complex domains $(B^n, J)$ and $(B^n, J')$ of this kind, which  admit $\cC^1$-extensions up to  the boundary (see also \cite{SS}).  Later, Gaussier and Sukhov proved  in \cite{GS} showed that the hypothesis of $\cC^1$-extendibility can be removed and that the result holds true for any pair of smoothly bounded, strictly pseudoconvex almost complex domains, proving Fefferman's theorem in almost complex setting in full generality (see also \cite{CGS2}).  
   \par
\smallskip
Motivated by these results and  possible applications on  plurisubharmonic exhaustions,  in this paper we determine more general situations in  which the existence of   foliations by stationary disks (and hence of generalized Riemann maps) is granted. Basically, we follow the approach of \cite{CGS}. We first   consider the differential problem that 
characterizes the stationary disks of  an almost complex domain $(D, J_o)$ and we explicitly determine  the  associated linearized operator $\gR$ at a given stationary disk $f_o: \overline \D \to \overline D$. When $\gR$ is invertible, we say that {\it $\partial D$ is  good relatively to the pair $(f_o, J_o)$\/}. A direct application of the  Implicit Function Theorem implies that if $\partial D$ is good, then there exists  stationary disks in a neighborhood of $f_o$  also when $J_o$ is replaced by a sufficiently close   almost complex structure $J$. 
On the base of this observation, one has  that {\it if an almost complex domain $(D, J_o)$  has a  foliation $\cF^{(x_o)}$ of stationary disks through $x_o$,  and if  the boundary is good for  $(f_o, J_o)$ for any $f_o \in \cF^{(x_o)}$,  then there exists  a  foliation  for $(D, J)$ of stationary disks passing through $x_o$,   also when  $J_o$ is replaced by  a  sufficiently close almost complex structure $J \neq J_o$\/}. \par
Secondly,  by  a line of arguments  that  goes back to Lempert and Pang  (\cite{Le, Pa}; see also \cite{Tr, ST, CGS, SS})),  we are able to prove  that  {\it any smoothly bounded, strictly linearly convex domain $D \subset \bC^n$ has a boundary which is ``good''  for any 
of its stationary disks\/}.  This fact and  previous observation  bring directly  to our result, which generalizes the quoted Coupet, Gaussier and Sukhov's  theorem on the unit ball:  {\it if a smoothly bounded, strongly pseudoconvex  domain  $D$ in an almost complex manifold  $(M,J)$  is  biholomorphic  to  a strictly linearly convex domain $\wh D \subset (\bC^n, J')$,  endowed with small deformation   $J'$   of  $\Jst$, then there exists  a  foliation by stationary disks of   $(D, x_o)$ for any $x_o \in D$\/} (Theorem \ref{main22}). \par
\smallskip
This shows that the class of almost complex domains,  admitting a   foliation by stationary disks,  is indeed much larger than the class  considered in \cite{CGS}. In fact, via  a diffeomorphism $\varphi: \cU \to \cV \subset \bC^n$ mapping $\overline{D}$ onto $\overline{B^n}$,  one obtains   the existence of  foliations by stationary disks on $(B^n, J' )$ also  when  $J' = \varphi_*(\wh J)$ is not   a small deformation of $\Jst$.\par 
We also prove  that,  for any almost complex domain $(D, J)$  as above  and with $J'$ sufficiently close to $\Jst$,  there exists a  generalized Riemann map $\varphi: \overline{B^n} \to \overline{D}$ for any  $x_o \in D$ and the  function $u= (\varphi^{-1})^*(u_o): D \to ]-\infty, 0[$ of   $u_o(z) \= \log(|z|)$  is a plurisubharmonic exhaustion  for $D$.  When  $J$ is integrable,   $u$ is a solution of the Monge-Amp\`ere equation $(\partial \bar \partial u)^n = 0$ with boundary data $u|_{\partial D} = 0$ and logarithmic singularity at $x_o$. It would be interesting to know if  this and other related properties  have counterparts in almost complex setting. 
\par
\smallskip
Finally,  we consider  the families $\cG^{(x_o, a)}$, formed by all  stationary disks in a given almost complex domain $(D, J)$, passing through  a given boundary point $x_o \in \partial D$ and with tangent vector $v$ at $x_o$, with inner product $<v, \nu> $  with the unit normal $\nu_{x_o}$  larger than a value $a \geq 0$. In case $D$ is a strictly convex domain in $\bC^n$, the disks in $\cG^{(x_o, a)}$ give a (regular) foliation of a certain subdomain $D^{(x_o, a)} \subset D$ that coincide with $D$ in case $a = 0$ (\cite{CHL}). We prove that  {\it if $a >0$, this is true also  when the standard complex structure $\Jst$  is replaced by an almost complex structure $J$ sufficiently close to $\Jst$\/} and we therefore have an analogue of the previous results also for what concerns  foliations of conical subdomains $D^{(x_o, a)}$, $a > 0$,  of almost complex domains.  A proof  for the  case $a = 0$ seems to be at the moment  out of reach,  because the family of stationary disks $\cG^{(x_o, 0)}$ is not parameterized by a compact set, in contrast with all other considered situations. \par
As final remark,  notice that when $J$ is integrable, the regular foliations $\cG^{(x_o, a)}$ determine analogues of the Riemann map and have been used in \cite{BP, BPT} to construct  solutions to 
the Monge-Amp\`ere equation $(\partial \bar \partial u)^n = 0$ with singularity at a given boundary point. 
It would be interesting to know if a similar construction can be obtained in an almost complex setting.\par
\smallskip
The structure of the paper is as follows. In \S 2 we recall a few basic facts  and the definition of stationary disks in almost complex domains. In \S 3,  we consider the so-called foliations of circular type,  prove their stability under small deformations of $J$   in case of  a ``good boundary''. In \S 4,  general conditions for a boundary ``to be good" are given and are used to show that any strictly linearly convex domain has a ``good''  boundary. This and the results of \S 3 give our main Theorem \ref{main22} as immediate consequence. Section \S 5 is devoted  to the quoted results on foliations of conical subdomains.\par
\bigskip
\section{Preliminaries}
\setcounter{equation}{0}
\subsection{Notations}
\label{notation}
Given a real manifold $M$ and a system of coordinates $\xi = (x^i): \cU \subset M\to \bR^n$, we call 
 {\it associated coordinates on $T^*M$ \/} the coordinates  $\hat \xi = (x^i,  p_i)$, where for any $\alpha  \in T_x^*M$ the ``$p_i$'' are the  components of $\a = p_i dx^i$ in the basis $(dx^i)$. If   $(M, J)$ is an almost complex manifold of real dimension $2n$, we  call   {\it system of complex coordinates\/} any local  diffeomorphism $\xi = (z^i): \cU \subset M \longrightarrow \bC^n$. 
We call them {\it holomorphic}  whenever $J$ is integrable and $\xi = (z^i)$ is a chart of   the corresponding complex manifold structure of $(M, J)$.  We also call {\it associated complex coordinates on $T^*M$} the complex  complex coordinates  $\wh \xi = (z^i,  w_j): \pi^{-1}(\cU) \subset T^*M \to \bC^{2n}$, where the $w_i$'s are defined for any 1-form $\a$ by the expression $\a = w_i dz^i + \overline{w_i} d\overline{z^i}$.\par
% If $J$ is an integrable complex structure and  $\xi = (z^i)$ maps $J$ into   the standard  complex structure $\Jst$ of $\bC^n$,   the complex coordinates $\xi = (z^i)$ are called    {\it holomorphic\/}.
 \par
%\smallskip
For any Banach space $X$ and   $\cU \subset \bR^M$,   
 $\a \in ]0 , 1[$, we denote by
$\cC^{\alpha}(\cU,X)$ the Banach space of the functions
$f: \cU \rightarrow X$ such that
$$ 
\parallel f \parallel_{\alpha} \=
\sup_{\zeta \in \cU} \parallel f(\zeta) \parallel +
\sup_{\theta,\eta \in \cU, \theta \neq \eta}
\frac{\parallel f(\theta) - f(\eta)\parallel}{\vert \theta - \eta
\vert^{\alpha}} < \infty.
$$
If $\alpha = m +
\beta$, for some $m \in \bN$ and $\beta \in ]0,1[$, we denote by  $\cC^{\alpha}(\cU, X)$ 
the Banach space
$ 
\cC^{\alpha}(\cU, X) =\{r \in \cC^m(\cU,X): D^{\nu} r \in C^{\beta}(\cU, X),
\nu: \vert \nu \vert \leq m\}
$. 
Finally, for any $\alpha, \epsilon>0$, we  set  
$\cC^{\alpha, \epsilon}(\bar 
\Delta, \bC^n) = \cC^\epsilon(\overline \Delta, \bC^n) \cap \cC^\alpha(\Delta, \bC^n)$ and 
$H^{ \epsilon}(\bar 
\Delta, \bC^n) = \cC^\epsilon(\overline \Delta, \bC^n) \cap Hol(\Delta, \bC^n)$
\par
\medskip
\subsection{ Lifts of $J$-holomorphic disks} We recall that a $\cC^\alpha$-map 
$f: M \to M'$, $1 \leq \a$, between two almost complex manifolds $(M,J)$, $(M',J')$  is called  {\it $(J,J')$-holomorphic\/} if  and only if $ \bar\partial_{J,J'}f (v) = 0$ for any $v \in TM$, 
where $ \bar \partial_{J,J'}f$ is the operator
\beq \label{d-bar} \bar \partial_{J,J'}f: TM \to TM'\ , \qquad \bar \partial_{J,J'}f (v) \= f_*\left(J(v) \right)- J'(f_*(v))\ .\eeq
When  $(M, J) = (\bC^n, \Jst)$, we  will  shortly  write  $\bar \partial_{J'}$ for $\bar \partial_{\Jst, J'}$. 
A  {\it $J$-holomorphic disk\/} of   $(M, J)$ is a   $(\Jst, J)$-holomorphic map $f: \Delta\to M$  from the unit disk $\Delta \subset \bC$ into $(M,J)$.  Recall that    $\bar \partial_J f = 0$   if and only if   $\bar \partial_J f \left(\left.\frac{\partial}{\partial x}\right|_{x + i y}\right) = 0$ at any $x + i y  \in \Delta$ (see e.g. \cite{IS}).  \par
\smallskip
If  $(M,J)$ is a  complex manifold,  the  cotangent bundle $T^*M$ 
 is naturally  endowed with   an  {\it integrable\/} complex structure $\bJ$, determined  by the identifications of  open subsets $\cU \subset M$ with  open subsets of $\bC^n$ and by the identifications of   the sets $T^*M|_{\cU}$ with open subsets of $\bC^{2n} = T^* \bC^n$.  When $J$  is  not integrable, these identifications are no longer valid, but 
 there still exists  a natural almost complex structure $\bJ$ on $T^*M$,   which reduces to the usual one  if  $J$ is integrable (\cite{IY}). The main properties  of $\bJ$ are summarized in the next 
proposition.  Here,  $J^i_j = J^i_j(x)$ are  the  components of  
$J = J^i_j \frac{\partial}{\partial x^i} \otimes dx^j$ in a  system of real coordinates $\xi = (x^i)$.\par
%\medskip
\begin{prop}\cite{IY} \label{liftedalmostcomplex} For any almost complex manifold $(M,J)$, 
there exists a unique almost complex structure $\bJ$ on $T^*M$ with the 
following properties:
\begin{itemize}
\item[i)] the projection $\pi: T^*M \to M$ is $(\bJ,J)$-holomorphic;
\item[ii)]  for any $(J,J')$-biholomorphism $f: M\to N$ between two almost 
complex manifolds $(M,J)$ and $(N,J')$, the induced map $\hat f: T^*N \to T^*M$
is $(\bJ',\bJ)$-holomorphic;
\item[iii)] if $J$ is integrable, then also $\bJ$ is integrable and coincides with 
the natural complex structure of $T^*M$;
\item[iv)]  in a  system 
of coordinates 
\beq \label{coordinates}
\wh \xi = (x^1, \dots, x^{2n}, p_1, \dots, p_{2n}): \pi^{-1}(\cU) \subset T^*M \longrightarrow \bR^{4n}\ ,\eeq   associated with  $\xi = (x^i)$, the tensor $\bJ$   is of the form 
$$ \bJ= J_{i}^a \frac{\partial}{\partial x^a} \otimes dx^i + 
J_{i}^a   \frac{\partial}{\partial p_i}\otimes d p_a  +  \phantom{aaaaaaaaaaaaaaaa}
 $$
\beq \label{bJ}  \phantom{aaaaaaaa} + 
\frac 1 2  p_a\left( - J^a_{i,j}  + J^a_{j,i} + J^a_\ell \left(J^\ell_{i,m} J^m_j - J^\ell_{j,m} J^m_i \right)\right)
\frac{\partial}{ \partial p_j}\otimes d x^i\ .\eeq
\end{itemize}
\end{prop}
The almost complex structure   $\bJ$  is called {\it  canonical 
lift of $J$ on $T^*M$\/}. \par
\begin{lem} \label{littlelemma} Let $\bJ$ be the canonical lift on $T^*M$ of an almost complex structure $J$.  For any $0 \neq t \in \bR$,   the map  $\varphi_t: T^*M \to T^*M$ defined by 
$\varphi_t(\alpha) = t \cdot \alpha$  is a $\bJ$-biholomorphic diffeomorphisms, i.e. 
 $\varphi_{t*} \circ \bJ = \bJ \circ \varphi_{t*}$.
\end{lem}
\begin{pf} Writing $\varphi_t$ in a system of coordinates   (\ref{coordinates}),  one has  that  $\varphi_t(x^i,  p_{j}) = (x^i,  t p_j)$. Using  (\ref{bJ}), the claim is then immediately checked.
\end{pf}
Given a $J$-holomorphic disk $f: \Delta \to (M, J)$, we call {\it lift of $f$\/} any $\bJ$-holomorphic disk $\wh f: \Delta \to (T^* M, \bJ)$ so that $f = \pi \circ \wh f$.\par
\medskip 
\subsection{Stationary disks} Let   $\Gamma \subset M$ be a smooth hypersurface  of an  almost complex manifold $(M,J)$. The {\it conormal bundle  of $\Gamma$\/} is  defined as 
\beq \cN \= \{\ \alpha \in T^*_xM\ ,\ x\in \G\ :\ \alpha|_{T_x\Gamma} \equiv 0\ \}  \subset T^* M|_{\G} \ .\eeq
In the following, we denote by  $\cN_* = \cN \setminus \{\text{zero section}\}$ and  {\it when we mention ``the conormal bundle''  we will always mean  $\cN_*$\/}.\par
The CR structure of  $\Gamma$ is  defined as the pair $(\cD, J)$ given by  the   distribution  \beq \cD = \bigcup_{x \in \Gamma} \cD_x \subset T\G\ ,\qquad \cD_x \= \{\ v\in T_x \G\ :\ J(v) \in T_x\G\ \} \eeq
endowed with  the family $J = \{J_x\}$ of complex structures $J_x \= J|_{\cD_x}$.  A {\it defining 1-form\/} for $\cD$    is a  1-form on $\G$ 
so that $\ker \vartheta|_x = \cD_x$ for any $x \in \G$. The  {\it Levi form at $x$\/}  is the quadratic form
$\cL_x: \cD_x \to \bR$ defined by $ \cL_x(v) \= - d\vartheta_x(v, Jv)$ for any $v\in \cD_x$ and (up a scalar  factor) it  is independent  on the choice of  $\vartheta$.   This last property  follows immediately  from the fact   that  for  any  vector field  $X^{(v)} \in \cD$ so that $X^{(v)}_x = v$ one has 
\beq \label{ecco} \cL_x(v) = - d\vartheta_x(X^{(v)}, JX^{(v)})   = \vartheta_x([X^{(v)}, JX^{(v)}])\ . \eeq
An {\it oriented\/}  hypersurface $\Gamma\subset M$ is called {\it strongly pseudoconvex\/}  if   $\cL_x$ is positive definite at every $x\in \Gamma$  when  determined by a    defining 1-form $\vartheta$  with   $\vartheta_x(J n) > 0$  for any   $n$
pointing in the  ``outwards'' direction. If $D \subset M$ is a bounded domain with smooth boundary $ \partial D$, we say that $D$ is {\it strongly pseudoconvex\/} when $\partial D$, oriented so that the  "outwards" directions are pointing outside $D$, is strongly pseudoconvex. \par
\smallskip 
 The following notion of ``stationary disk''  for domains in  almost complex manifolds  was considered for the first time  by Coupet, Gaussier and Sukhov in \cite{CGS}. It generalizes the notion of  stationary disks of bounded domains in $\bC^n$ (\cite{Le, Tu}).\par
%\smallskip
\begin{definition} \label{stationarydisks} 
Let $D \subset M$ be a domain  with smooth boundary and $\cN_*$ the conormal bundle of $\partial D$. Given $\a\geq 1$, $\varepsilon> 0$, a  map  $f: \bar \Delta \to M$ is called {\it $\cC^{\a,\varepsilon}$-stationary disk  of $D$\/} 
if 
\begin{itemize}
\item[i)] $f|_{\Delta}$ is  a $J$-holomorphic embedding and  $f(\partial \Delta) \subset \partial D$; 
\item[ii)]  there exists a  lift $\wh f :  \overline \Delta \to T^*M$ of $f$  so that 
\beq \label{stationarycondition} \zeta^{-1}\cdot \hat f(\zeta) \in \cN_*\ \  \text{for any} \ 
\zeta \in \partial \Delta\eeq
 and  $ \wh \xi \circ \wh f \in \cC^{\a, \varepsilon}(\overline \D, \bC^{2n})$ for some complex coordinates $\wh \xi = (z^i, w_j)$  around  $\wh f(\overline \D)$. 
Here  ``\ $ \cdot$\ "  denotes   
the usual $\bC$-action  on $T^*M$, i.e. 
\beq\label{product} \z\cdot \alpha \= \Re(\z)\alpha  -  \Im(\z) J^*\alpha\quadÊ\text{for any}\ \alpha \in T^* M, \ \ \z \in \bC\ \ .\eeq
 \end{itemize}
\end{definition}
\medskip
In the following,  the values of $\a$ and $\varepsilon$ are considered as  fixed and by  ``stationary''  we always mean   ``$\cC^{\a,\varepsilon}$-stationary''. Moreover, 
given a  stationary disk $f$,  the maps  $\hat f$ satisfying   (ii)
are  called {\it  stationary lifts of $f$\/}.\par
%\medskip
%Here are two  direct properties of $J$-holomorphic disks and stationary disks.  \par
\begin{lem} \label{lemmino}  
i)  If $D \subset M$ is  a smoothly bounded,  strongly pseudoconvex  domain and $f: \overline\D \to M$ is a non-constant stationary disk of $D$, then $f(\overline{\D})  \subset \overline{D}$ and $f(\zeta) \in \partial D$ if and only if $\zeta \in \partial \D$.\par
\noindent ii)  For any $t \in \bR_*$ and any stationary lift  $\wh f$ of a  stationary disk  $f: \overline \D \to \overline D$, also the map $\wh f_t(\zeta) \= (\varphi_t \circ \wh f)(\z) = t \cdot \wh f(\zeta)$ is a stationary lift of $f$.
\end{lem}
\begin{pf} (i) If $D$ is strongly pseudoconvex, it is known that there exists  a defining function $\rho: \cU \subset M \to \bR$  for $D$ which is $J$-plurisubharmonic,   i.e.   so that 
$\rho \circ f: \D \to \bR$ is strictly subharmonic  for any $J$-holomorphic disk $f: \D \to \cU$   (see e.g. \cite{CGS2}, p.14). Since $\rho \circ f|_{\partial \D} =  0$,  the claim follows  from the maximum principle. \par
(ii) It  follows from the fact  that $\wh f_t$ satisfies (\ref{stationarycondition}) and that the diffeomorphism $\varphi_t$  is a   $\bJ$-biholomorphism  by Lemma \ref{littlelemma}.
\end{pf}
We conclude recalling the following theorem that   generalizes a well-known result by  Webster to the almost complex setting  (\cite{Web}). \par
%\medskip
\begin{theo} \label{spiro} \cite{Sp} Let $\Gamma$ be a strongly pseudoconvex hypersurface in 
an almost complex manifold $(M,J)$ and $\cN_* \subset T^*M$ its conormal bundle with the zero section excluded. Then  $\cN_* $  is a totally real submanifold of $(T^*M, \bJ)$.
\end{theo}
 \bigskip
\section{Foliations by stationary disks\\
and deformations of almost complex structures}\par
\setcounter{equation}{0}
\subsection{The Riemann-Hilbert problem for  stationary disks}In this and the next sections, $D$ is a  strongly pseudoconvex domain in  an almost complex manifold $(M, J)$ with smooth boundary $\partial D$ with   conormal bundle  $\cN \subset T^*M|_{\partial D}$.
We also assume that $\overline D \subset M$ is contained in a globally coordinatizable open subset $\cU \subset M$ or, equivalently, that $D$ is a domain of  $M = \bR^{2n} \simeq \bC^n$ equipped with a non-standard complex structure $J$.  We also assume that $D$ has  a smooth defining function $\rho: \cU \subset M \to \bR$ on $\cU$, so that 
$$D = \{\ x \in M \ :\ \rho(x) < 0\ \}\qquad \text{and}\qquad d \rho_x \neq 0\qquad\text{for any} \ x \in \G = \partial D\ .$$
\smallskip
We want to study the differential problem that characterizes the lifts   $\wh{f}: \overline \D \to T^*M$ of stationary disks  of $D$. First of all, consider 
the map 
\beq \wt \rho:  \bR_* \times T^*M|_{\cU}  \longrightarrow \bR \times T^*M|_{\cU}\ ,\ \  \label{defN} \wt \rho(t, \a) \= (\rho(\pi(\a)), \a - t \cdot d\rho_{\pi(\a)})\ .\eeq
Notice   that the bundle $\cN_* = \cN \setminus \{\text{zero section}\}$, which is a $2n$-dimensional submanifold of $T^*M$,  can be identified with the level set $$\{(t, \a):\ t \neq 0\ , \  \wt \rho(t, \a) = (0_\bR,0_{{\phantom{a}_{\!\!\!\! T^*_{\pi(\a)}M}}})\} \subset \bR_* \times T^*M|_{\cU}\ ,$$ 
which is a $2n$-dimensional submanifold of $\bR_* \times T^*M$. Therefore, using  a system of coordinates $\wh \xi = (x^i,p_j)$ on $T^*M|_{\cU}$, associated with  coordinates $\xi = (x^i)$,  we may  identify  $\bR_* \times T^*M|_{\cU}$ with an open subset $\cV \subset \bR^{4n+1}$  and   $\cN_*$  with the level set in $\cV$ defined by 
$$\cN_* \simeq \{\ (t, \a) \in \cV\ :\ \wt \rho^i(t, \a) = 0\ ,\ \ 1 \leq i  \leq 2n + 1\}\ .$$
By a direct check of the rank of the Jacobian, one can check that the  map $\wt \rho = (\wt \rho^1, \dots, \wt \rho^{2n+1})$
is  a smooth   defining function    for $\cN_*$.\par
\smallskip
We now  consider the map 
$ÿr: \bC \times\cV\subset \bC \times \bR^{4n+1} \longrightarrow \bR^{2n+1}$, defined by 
\beq \label{definitionr}  r(\z, t, \alpha) \= \left(\wt \rho^1(t, \zeta^{-1} \cdot \alpha), 
\dots, \wt \rho^n(t, \zeta^{-1} \cdot \alpha)\right)\ .\eeq
Here,  the product $\zeta^{-1} \cdot \alpha$ is as in  \eqref{product}.
By definition,  a disk  $f: \overline \D \to \overline D \subset \bR^{2n}$ is stationary if and only if there exists    $\wh f \in \mathcal (\cC^{\alpha, \epsilon}(\bar{\Delta}); \bC^{2n})$  and  $\lambda \in \cC^{\epsilon}(\partial \Delta; \bR)$ so that
\smallskip
\beq \label{Riemann-Hilbert}
\left\{
\begin{array}{lll}
\overline{\partial}_{\bJ}\wh f(\zeta) &=& 0, \ \ \ \zeta \in \Delta\\
& &\\
r(\zeta,  \lambda(\zeta), \wh f(\zeta)) &=& 0, \ \ \ \zeta \in \partial \Delta
\end{array}
\right.
\eeq
where $\overline\partial_\bJ = \overline\partial_{\Jst, \bJ} :
(C^{\alpha}(\overline\Delta); \bC^{2n}) \longrightarrow (C^{\alpha -
  1}(\overline\Delta; \bC^{2n}))$ is the operator   (\ref{d-bar}).\par
The differential problem (\ref{Riemann-Hilbert}) belongs to a class   often called  {\it of generalized Riemann-Hilbert problems} (see f.i.  \cite{MP}, Ch. VII).\par
%\bigskip
\subsection{Stability under small deformations of the data}
\label{stability}
Consider  a fixed almost complex structure   $J = J_o$,  a point $x_o \in D (\subset \bR^{2n})$ and a vector  $v_o \in T_{x_o} D \simeq \bR^{2n}$ and  denote by $\cR_{(J_o, x_o, v_o)} = (\cR_1, \dots, \cR_5) $ the operator from
$\cC^{\alpha, \varepsilon}(\overline\Delta; \bC^{2n})  \times  \cC^{\epsilon}(\partial \Delta; \bR) \times \bR_* $
into $ 
\cC^{\alpha-1, \varepsilon}(\overline\Delta; \bC^{2n})  \times   \cC^{\epsilon}(\partial \Delta; \bR^{2n+1}) \times \bC^{n}\times \bC^n \times \bR$ with components $\cR_i$ 
defined by
$$\cR_1(\wh f, \l, \mu ) \= \overline\partial_{\bJ_o}\wh f\ \ ,\ \ 
\cR_2(\wh f, \l, \mu ) \= r(\zeta,  \lambda(\zeta), \wh f(\zeta))\ ,$$
$$  \cR_3(\wh f, \l, \mu ) \=\pi(\wh f)|_{\zeta = 0} - x_o\ ,\ \ \cR_4(\wh f, \l, \mu ) \= \pi(\wh f)_*\left(\left.\frac{\partial}{\partial x} \right|_{\zeta = 0}\right) - \mu v_o\ ,$$
\beq \label{nonlinearoperator} \cR_5(\wh f, \l, \mu ) \=  \wh f \left(\pi(\wh f)_*\left(\left.\frac{\partial}{\partial x} \right|_{1}\right)\right) - 1\ .\eeq
Notice also that,  by Hopf's Lemma and Lemma \ref{lemmino} (ii), for any stationary disk, there exists a stationary lift  satisfying   $\wh f \left(\pi(\wh f)_*\left(\left.\frac{\partial}{\partial x} \right|_{1}\right)\right) = 1$. So, by   the previous section, the existence of a stationary disk $f: \D \to D$ with $f(0) = x_o$ and $f_*\left(\left.\frac{\partial}{\partial x} \right|_0\right) \in \bR v_o$ is equivalent to the existence of a solution  to\beq \label{theproblem} \cR_{(J_o, x_o, v_o)}(\wh f, \l, \mu) = 0\ .\eeq
Let $(\wh f_o, \l_o, \mu_o)$ be solution of (\ref{theproblem}) and  
$\gR_{(J_o, x_o, v_o; \wh f_o, \l_o, \mu_o)} \=\dot  \cR_{(J_o, x_o, v_o)}|_{(\wh f_o, \l_o, \mu_o)}$
 the linearized operator at $(\wh f_o, \l_o, \mu_o)$ determined by $ \cR_{(J_o, x_o, v_o)}$.
Now, by  the  Implicit Function Theorem (see e.g. \cite{KA}), when   $\gR = \gR_{(J_o, x_o, v_o; \wh f_o, \l_o, \mu_o)}$ is  invertible,   there exists a solution to the problem $\cR_{(J_t, x_t, v_t)}(\wh f, \l, \mu) = 0$ 
 for any smooth deformation $(J_t, x_t, v_t)$ of $(J_o, x_o, v_o)$ for $t$ sufficiently small $t$ and  $\dim_\bR \ker \gR_{(J_o, x_o, v_o; \wh f_o, \l_o, \mu_o)}$ is equal to 
 the dimension of the solutions space.  This motivates the following:\par
 \begin{definition} 
 Let $f_o: \overline\D \to \overline D$ be  a stationary disk of $(D, J_o)$   with $x_o = f(0)$ and $v_o = (\wh f)_*\left(\left.\frac{\partial}{\partial x} \right|_{\zeta = 0}\right)$. 
 We  call  $\partial D$ {\it  a good boundary for  $(J_o, f_o)$\/}  if there is a lift $\wh f_o$ of $f_o$ and a function $\l_o$ so that    $(\wh f_o, \l_o, 1)$ is a solution to (\ref{theproblem}) and the  linearized operator $\gR = \gR_{(J_o, x_o, v_o; \wh f_o, \l_o, 1)}$ is invertible.
 \end{definition}
%\bigskip
The Implicit Function Theorem and previous remarks brings immediately to the next proposition.  In the statement,  we denote by $g$  a  fixed Riemannian metric $g = g_{ij} dx^i \otimes dx^j$ on a neighborhood  of $\overline D$ and  by $g^* = g_{ij} dx^i \otimes dx^j + g^{ij} dp_i \otimes dp_j$  the corresponding Riemannian metric on  $T^*M$. We also set
\beq \label{norm} \| J - J'\|^{(1)}_{ \overline D}\= \sup_{x \in \overline D, v \in T(T^*_x M)} \frac{\| \bJ(v) - \bJ'(v)\|_{g^*}}{\|v\|_{g^*}}\ ,\eeq
where $\| \cdot \|_{g^*}$ is the norm function determined by $g^*$.  The topology determined  by  the norm $\| \cdot \|_{ \overline D}^{(1)}$ is clearly independent on the choice of $g$.  \par
%\medskip
\begin{prop} \label{main3}
 Let $f_o: \overline \D \to \overline D$ be a stationary disk of $D \subset (M, J_o)$  with $x_o = f_o(0)$ and $v_o=f_o{}_*\left(\left.\frac{\partial}{\partial x} \right|_{\zeta = 0}\right) $.
 If $\partial D$ is a good boundary  for  $(J_o, f_o)$,   there exists a neighborhood $\cV \subset D$ of $x_o$,  a neighborhood $\cW \subset T D$ of $v_o$, with $\pi(\cW) = \cV \subset D$
 and a real number  $\varepsilon > 0$ so that,  for any $x \in \cV$, $v\in \cW$  and    $\| J - J_o\|^{(1)}_{\overline D} < \varepsilon$,  there exists a unique stationary disk $f$ of $(D, J)$  so that 
\beq \label{pinco} f(0) = x\ ,\qquad f_*\left(\left.\frac{\partial}{\partial x}\right|_{\zeta = 0}\right) = \mu v\ \ \ \text{for some}\ \ \mu \neq 0\ .\eeq
The disk $f$ depends differentially on  $x$,   $v$ and   $J$  and,  given $m_o > 0$,  one can choose $\varepsilon$, $\cW$ and $\cV = \pi(\cW)$ so that  $\sup_{\zeta \in \overline \D} \dist_g(f(\zeta), f_o(\zeta))  < m_o$.
\end{prop}
\subsection{Foliations of circular type and their stability}
\subsubsection{Blow-up of an almost complex domain at one point}
Let $x_o$ be a point of  the almost complex manifold $(M, J)$ and $\xi = (z^i): \cU \to \bC^n$  a system of complex coordinates  with
 \beq
  \label{coord} 
 \xi(x_o) = 0 \ ,\qquad \xi_*(J|_{x_o}) = \Jst|_0 \ . \eeq 
Consider the blow up $\pi: \wt \cU \to \xi(\cU)  \subset \bC^n$ of $\xi(\cU) $ at $0$, i.e. the submanifold of $\bC^n \times \bC P^{n-1}$ defined by $\wt \cU = \{\ (z, [w])\ : z \in [w]\ , z \in \cU\}  \subset \bC^n \times \bC P^{n-1}$. The standard projection $\pi(z, [w]) = z$  composed with $\xi^{-1}$ determines  a diffeomorphism between  $\wt \cU \setminus \pi^{-1}(0)$ and $\cU \setminus \{0\}$ that we  use  to glue $\wt \cU$ with $M \setminus \{x_o\}$ and obtain a manifold $\wt M$ that we call {\it  blow up of $(M, J)$ at  $x_o$\/}. \par
%\smallskip
At a first glance, this  construction  seems to depend on the  choice of the complex coordinates $\xi = (z^i)$.   But indeed the 
 real manifold structure of $\wt M$  depends only  on the linear map  $J_{x_o}: T_{x_o} M \to T_{x_o} M$. This fact is a direct consequence of the following simple lemma.
\begin{lem}  Consider two sets of  complex  coordinates  $\xi = (z^i)$ and $\xi' = (z'{}^j)$  on $\cU$ satisfying 
$\xi(x_o) = \xi'(x_o) = 0$ and $\xi_*(J_{x_o}) = \xi'_*(J_{x_o}) = \Jst|_0$. Then the diffeomorphism  $\wt \varphi = \pi^{-1} \circ (\xi' \circ \xi^{-1}) \circ \pi$  of  $\wt \cU \setminus \pi^{-1}(0)$ into itself 
%defined by 
%$
%:  \wt \cU \setminus \pi^{-1}(0) \longrightarrow  \wt \cU \setminus \pi^{-1}(0)\ ,
%$$ 
admits a unique smooth extension on  $\wt \cU$.  It follows that the blow up $\wt M$,  defined using the chart $\xi = (z^i)$,  is naturally diffeomorphic to the one constructed using the chart $\xi' = (z'{}^i)$. 
\end{lem}
\begin{pf}
%The claim on $\wt \varphi$ can be checked as follows. 
By construction, the map  $\varphi = \xi'\circ \xi^{-1}$ is so that $\varphi_*|_0 \circ  \Jst = \Jst  \circ \varphi_*|_0$ and hence it  is of the form 
\beq \label{approximation}  \varphi(z) = \psi(z) +  g(z)\eeq
where   $\psi$ is the $\bC$-linear map $\psi = \varphi_*|_0:   \bC^n \to \bC^n$ and $g: \cU \to \cU$
is  an infinitesimal of the second order in $|z|$.  Since  
$$\wt \varphi(z, [z])  =  (\pi^{-1} \circ \varphi \circ \pi)(z, [z])  = (\psi(z) + g(z), [\psi(z) + g(z)])\ ,$$ 
an explicit computation  in coordinates shows that  $\wt \varphi$ extends smoothly on  $ \pi^{-1}(0) \subset \wt \cU$ by setting $\wt \varphi(0, [v])  \= (0, [\psi(v)])$ for any $[v] \in \bC P^{n-1}$.
\end{pf}
\smallskip
\subsubsection{Foliations of circular type}
Let $D$ be a smoothly bounded, strongly pseudoconvex domain  in  $(M, J)$ and denote also by $\overline{\wt D} \subset \wt M$   the  blow up of  $\overline{D}$ at a point $x_o$ as defined in the previous section. 
For any  stationary disk $f: \overline{\D} \to \overline{D}$ with $f(0) = x_o$ and  $
f_*\left( \left.\frac{\partial}{\partial x}\right|_0 \right) = v$ there exists a unique map $\wt f:\overline{\D}\to\overline{ \wt D}$ so that $\pi \circ f(\zeta) = f(\zeta)$ for any $\zeta \neq 0$. In fact, if we identify $\wt D$ with  a domain in $\wt \cU \subset \bC^n \times \bC P^{n-1}$ by means of  a chart  like in (\ref{coord}),  the lifted map $\wt f$ is of the form
\beq \wt f(\zeta) = \left\{\begin{array}{ll} (f(\zeta), [f(\zeta)]) \ & \text{when} \ \zeta \neq 0\ ,\cr \phantom{a}\cr  (0, [v])\ & \text{when} \ \zeta = 0\ .\end{array}\right.\eeq
 Since $f$ is $J$-holomorphic (and hence  $f_*(\Jst|_0) = J|_0 = \Jst|_0$),  we may write
 \beq f(\zeta) = h(\zeta) + g(\zeta)\eeq
for some  holomorphic disk $h: \overline{\D} \to \cU \subset \bC^n$ and a smooth map $g: \overline{\D}\to \cU$ which is  infinitesimal of second order in $|\zeta|$. Using this,  one can check that $\wt f$ is smooth also at $0$.  We call $\wt f$ {\it the smooth lift of $f$ at  $\wt D$\/}.\par
\begin{definition} \label{foliation} Let   $x_o \in D$ and $\wt D$ as above and 
 denote by $\cF^{(x_o)}$ the family  of all stationary disks of $D$ with  $f(0) = x_o$. We call  $\cF^{(x_o)}$  {\it  foliation of circular type of the pointed domain $(D, x_o)$\/}  if  the following conditions are satisfied:
\begin{itemize}
\item[i)] for any $v \in T_{x_o} D$, 
 there exists a unique disk $f^{(v)} \in \cF^{(x_o)}$
 such that 
$f^{(v)}_*\left( \left.\frac{\partial}{\partial x}\right|_0 \right) = \mu \cdot v$ for some $0 \neq  \mu \in \bR$;
\item[ii)]  previous an identification  $(T_{x_o} D, J_{x_o}) \simeq (\bC^n, \Jst)$,  
 the map
\beq \label{exp} \exp: \wt B^n  \subset \wt \bC^n \longrightarrow  \wt D \ \ , \qquad  \ \exp(v, [v]) \= \wt{f^{(v)}}(|v|)\ ,\eeq
between the blow up  at $0$ of   $B^n\subset \bC^n $  and the blow up of $D$ at $x_o$ 
is smooth  with a smooth extension  up to the boundary, which induces a diffeomorphism  between the boundaries  $\exp|_{\partial B^n}: \partial B^n \to \partial D$.
 \end{itemize}
If  $\cF^{(x_o)}$ is  a  foliation of circular type,  we call   $x_o$   {\it center of the foliation\/}  and  $D$  a  {\it  domain  of  circular type w.r.t. to $J$}. \end{definition}
\bigskip
\subsubsection{Stability under small deformations of foliations of circular type }
\begin{prop} \label{main1} 
Let $D$ be of 
circular type w.r.t. to   $J_o$ and with center $x_o$. If   $\partial D$ is a good boundary for $(J_o, f_o)$ for  any stationary disk  $f_o \in \cF^{(x_o)}$, then there exists  $\varepsilon > 0$ and an open neighborhood $\cU \subset D$ of $x_o$  so that  for any $J$ with   $\| J - J_o\|^{(1)}_{\overline D} < \varepsilon$ and any $x \in \cU$, the point $x$  is  center  of a foliation of  circular type of $D$ w.r.t.  the almost complex structure $J$.
\end{prop}
\begin{pf} Using a system of coordinates $\xi = (x^i)$ on a neighborhood $\cW$ of $x_o$, let us  identify $\cW$ with an open subset of $\bR^{2n} \simeq \bC^n$ and its tangent space  with    $T \cW \simeq \cW \times \bR^{2n} \subset \bR^{4n}$.  Pick also the same
 Euclidean inner product $<, >$ on  all tangent spaces in  $T\cW \simeq \cW \times \bR^{2n}$. 
By  definitions, for any  $v_o \in S^{2n-1}_{x_o} = \{ \ v \in T_{x_o} M\ : \ < v,v> = 1 \ \}$, 
 there is a unique stationary disk $f \in \cF^{(x_o)}$   with $f_*\left(\left. \frac{\partial}{\partial x} \right|_0\right) = \mu \cdot v_o$ for some $\mu \neq 0$. \par
By  Proposition \ref{main3}, there exists a neighborhood 
$\cU^{(v_o)}$ of $x_o$,  a neighborhood $\cV^{(v_o)} \subset S^{2n-1}$ and $\varepsilon^{(v_o)} > 0$, so that, for any 
$y \in \cU^{(v_o)}$, $v \in \cV^{(v_o)}ÿ\subset T_y M \simeq  T_{x_o} M \simeq \bR^{2n}$ and $J$  with $\| J - J_o\|^{(1)}_{\overline{D}} < \varepsilon^{(v_o)}$, 
there exists a unique disk $\wt f$, which is stationary  for $D$ w.r.t. $J$, 
passing through $y$ and  with $\wt f_*\left(\left. \frac{\partial}{\partial x} \right|_0\right) $ parallel to $v$.
By compactness of $S^{2n-1}$,   there exists  a finite number of vectors $v_1, \dots, v_N \in S^{2n-1}$ so that the corresponding 
open sets $\cV^{(v_i)} \subset S^{2n-1}$ give an open covering of $S^{2n-1}$. We conclude that, for any point $y \in \wt \cU = \bigcap_{i=1}^N \cU^{(v_i)}$,   $\| J - J_o\|^{(1)}_{\overline{D}} <\min_i \varepsilon^{(v_i)}$ and  $v \in T_y M$,  there exists 
a unique disk passing through $y$, which is stationary w.r.t. $J$  and with  $\wt f_*\left(\left. \frac{\partial}{\partial x} \right|_0\right) $ parallel to $\frac{v}{|v|} \in S^{2n-1}$. In particular,  the  disks in $\cF^{(y)}$, $y \in \wt \cU$,  satisfy Definition \ref{foliation} (i). \par
Consider now the map $\exp: \overline{\wt B^n} \to \overline{\wt D}$ in  (\ref{exp}). By Proposition \ref{main3},  it  is smooth and depends smoothly on $y$ and $J$. Moreover, if $J = J_o$ and $y = x_o$,  it is a diffeomorphism between  manifolds with boundaries. Hence,   there exists  $\cU \subset \wt \cU$ and  $\varepsilon < \min_i \varepsilon^{(v_i)}$ so that   $\exp_*$ is invertible at all points of $\overline{\wt B^n}$  whenever $y \in \cU$ and   $\| J - J_o\|^{(1)}_{\overline{D}} <\varepsilon$. In these cases,  $\exp$ is  a local homeomorphism  from the compact set 
$ \overline{\wt B^n}$ to $\overline{\wt D}$ and hence  is a covering map of $\overline{\wt D}$. 
%(see e.g. \cite{DoC}, p. 374)
Being $\overline{\wt B^n}$ simply connected,  it is a diffeomorphism, i.e.  also (ii) of Definition \ref{foliation} holds true.
\end{pf}
\bigskip
\section{Conditions that force a boundary to be ``good''}
\label{goodness}
\setcounter{equation}{0}
In this section we are going to prove a result (Theorem \ref{main2}),  which provides a condition for the existence of foliations by stationary disks of a pointed domain $(D, x_o)$  endowed with a small deformation of the standard complex structure.  An immediate consequence of this and of the contents of \S 3 is  represented by the following theorem. \par
\begin{theo} \label{main22} Let $D \subset M$ be a  smoothly bounded, strongly  pseudoconvex 
domain in an almost complex manifold $(M, J_o)$. If there is a local diffeomorphism $\varphi: \cU \subset M \to \bC^n$, so that  $\wh D = \varphi(D)$ is  a strictly linearly convex domain $\wh D \subset \bC^n$  and $\varphi_*J_o$ is sufficiently close to $\Jst$ in $\cC^1$-norm, then $D$ is a domain of circular type w.r.t. $J$ and any point is a center. 
\end{theo}
Roughly speaking, this shows that  if one defines a suitable topology on the set of almost complex domains admitting foliations of circular type,  such space contains a whole open neighborhood of the class of  strictly linearly convex domains of $\bC^n$.\par 
\medskip

\subsection{The linearized operator $\gR = \gR_{(J_o, x_o, v_o; \wh f_o, \l_o, \mu_o)}$}\label{condcoer}
First of all, we want to determine an explicit expression for the tangent map $\gR = (\gR_1, \gR_2, \gR_3,\gR_4, \gR_5)$  at  $(\wh f_o, \l_o, \mu_o)$ of the operator  (\ref{nonlinearoperator}). For this, recall that, being  $\wh f_o: \overline \D \to T^*M$  a $\bJ$-holomorphic disk, one can always find  a  system of complex coordinates $(z^i)$ on  a neighborhood $\cW$ of $\wh f_o(\overline \D)$,   in such a way that, identifying  $\cW$ with an open subset  of $\bC^{2n}$, one has  $\bJ|_z = \Jst|_z$ at any $z \in \overline{f(\D)}$.  Moreover,   by Hopf lemma and being  the defining function $\rho$ strongly plurisubharmonic,  we have that $d \rho\left(f_o{}_*\left(x \frac{\partial}{\partial x} +y \frac{\partial}{\partial y} \right)\right) = 
d \rho\left(\Re\left(z^1 \frac{\partial}{\partial z^1}\right)\right) \neq 0$ at all points of $f(\partial \D)$.
In these coordinates, the tangent map of $\cR_1 = \overline\partial_{\bJ_o}$ at $\wh f_o$   is \beq
%\label{3.25}
 \gR_1(\wh h) = \frac{\partial \wh h}{\partial \overline\zeta} + 
\frac{1}{2i}D(\bJ - \Jst)_{\hat f_o} \cdot \wh h\ ,\eeq
where  $D(\bJ - \Jst)_{\hat f_o}$ is  the real differential of the matrix
valued function $\zeta \mapsto (\bJ - \Jst)_{\hat f_o(\zeta)}$.  In 
 matrix notation,  $D(\bJ - \Jst)_{\hat f_o} \cdot \wh h$ can be written as
 $$\left(D(\bJ - \Jst)_{\hat f_o}\cdot \wh h\right)_\zeta = A(\zeta) \cdot \wh h(\zeta) + B(\zeta)\cdot \overline {\wh h}(\zeta)\ ,$$
for some 
 $A, B: \overline \D \to M_{n \times n}(\bC)$ and  $ \gR_1$ assumes   the form
\beq 
 \gR_1(\wh h)  =
 \frac{\partial \wh h}{\partial \overline\zeta} + A \cdot \wh h + B\cdot \overline {\wh h}\ .\eeq
\par
% \medskip
 Consider now the tangent map  $\gR_2$.  By previous remarks,   the defining function $\wt \rho = (\wt \rho^1, \dots, \wt \rho^{2n+1})$  in    (\ref{defN})  is  locally equivalent to 
 $$\wh \varrho(t, \a) = \left( \varrho(\a), t - 
 \frac{\a\left(\Re\left(z^1\frac{\partial}{\partial z^1}\right)\right)}{\left.d \rho\left(\Re\left(z^1\frac{\partial}{\partial z^1}\right)\right)\right|_{\pi(\a)}}\right)$$
 where   $  \varrho: \cW \subset T^*M|_{\cU} \to \bR^{2n}$ is  the defining function for $\cN_*$ obtained  by   replacing $t =\frac{\a\left(\Re\left(z^1\frac{\partial}{\partial z^1}\right)\right)}{\left.d \rho\left(\Re\left(z^1\frac{\partial}{\partial z^1}\right)\right)\right|_{\pi(\a)}}$ in all places of $\wt \rho$.  
 If we set 
 \beq \label{varrhozeta} \wh r(\zeta, t, \a) \= \left(\varrho(\z, \a), t - \frac{\a\left(\Re\left(z^1\frac{\partial}{\partial z^1}\right)\right)}{\left.d \rho\left(\Re\left(z^1\frac{\partial}{\partial z^1}\right)\right)\right|_{\pi(\a)}}\right),\ \text{with}\  \varrho(\z, \a)\= 
\varrho(\z^{-1} \cdot \a)\eeq
we see that   $\gR_2$  is equivalent to the tangent map of the operator 
 $$\cR_2(\left.\wh f\right|_{\partial \D}, \lambda)=  \left( \varrho(\cdot,   \wh f(\cdot)), 
 \lambda - \frac{\wh f(\cdot)\left(\Re\left(z^1\frac{\partial}{\partial z^1}\right)\right)}{\left.d \rho\left(\Re\left(z^1\frac{\partial}{\partial z^1}\right)\right)\right|_{\pi(\wh f(\cdot))}}
 \right)$$
and hence  of the form 
\beq \gR_2(\wh h, \tau) =  \left(2 \Re (G \cdot \wh h|_{\partial \D}), \tau - g(\wh h)\right)\eeq
where $g$ is  obtained by linearization of the map
$\wh f \mapsto \frac{\wh f(\cdot)\left(\Re\left(z^1\frac{\partial}{\partial z^1}\right)\right)}{\left.d \rho\left(\Re\left(z^1\frac{\partial}{\partial z^1}\right)\right)\right|_{\pi(\wh f(\cdot))}}$ and 
 $G$ is the matrix valued map on $\partial \D$ defined by 
\beq \label{matrixG}
G(\zeta) = \left(
\begin{matrix}
\frac{\partial \varrho^1}{\partial  z^1}(\zeta,  \hat f(\zeta)) &\cdots&\frac{\partial
\varrho^1}{\partial z^{2n}}(\zeta,  \hat f(\zeta))\\
\vdots&\ddots&\vdots\\
\frac{\partial \varrho^{2n}}{\partial z^1}(\zeta,  \hat f(\zeta))& \cdots&\frac{\partial \varrho^{2n}}
{\partial  z^{2n}}(\zeta,  \hat f(\zeta))
\end{matrix}
\right)\ ,\quad  \z \in  \partial \Delta\eeq
By Theorem \ref{spiro},    $\cN_*$ is totally real  w.r.t.  $\bJ$ and 
hence, by  our choice of the coordinates,  it is  totally real also w.r.t. $\Jst$ on a neighborhood of $\hat f(\partial \D)$. This implies that 
\beq \label{Lopatinsky}
\det \left(G(\zeta)\right) \neq 0, \qquad \text{for any}\ \zeta \in \partial\Delta\ .
\eeq
Finally, the maps $\gR_3$, $\gR_4$ and $\gR_5$ are easily seen to be  (here $h \= \pi \circ \wh h$)
$$\gR_3(\wh h) = h(0)\ ,\ \gR_4(\wh h, \sigma) = \left.\frac{\partial  h}{\partial x}\right|_{\zeta = 0} - \sigma v_o\ ,\ \gR_5(\wh h) = \wh f_o\left(\left.\frac{\partial h}{\partial x}\right|_1\right) + \wh h
\left(\left.\frac{\partial f_o}{\partial x}\right|_1\right).$$\par
%\bigskip
 \subsection{The operator  $R_{A,B,G} \ = (\gR_1, \gR_2)$}
 \label{R_{A,B,G}}
Consider the operator 
$$R_{A,B,G}  = (\gR_1, \gR_2) = \left(\frac{\partial \wh h}{\partial \overline\zeta} + A \cdot \wh h + B\cdot \overline {\wh h}\ ,\  2 \Re (G \cdot \wh h)\right)\ ,$$
which  is a well-known Fredholm operator related  with the generalized
Riemann-Hilbert problems.  In the next theorem, we recall some information  that will be used in the sequel (see e.g. Thm. 3.2.5, Thm. 3.3.1 in \cite{Wen}).\par
%\medskip
\begin{theo} \label{indexR} If $G$ satisfies    (\ref{Lopatinsky}), the operator $R_{A,B,G}$ is Fredholm  with index 
$ \nu = 2n - \frac{1}{i \pi} \int_{\partial \Delta} d\arg(\det(G))$
and hence is surjective if and only if 
\beq\label{surjectivity} \dim \ker R_{A,B,G} = 
2n - \frac{1}{i \pi} \int_{\partial \Delta} d\arg(\det(G))\ .
\eeq
\end{theo}
%\bigskip
Next, we need to recall  a lemma due to Globevnik and some of its direct consequences, 
which give a way to establish the surjectivity of  $R_{A,B,G} $  in case of integrable complex structures.
But in order to state them,  we first need to recall the definition of ``canonical system'' (see e.g. \cite{Gl}). In what follows, for any holomorphic function $g: \cU \subset \bC \to \bC^N$
on a neighborhood of $\infty$ and with  at most one pole at $\infty$, we call {\it order of (zero of) $g$\/} the integer $k$ 
such that $g = \frac{1}{z^k} g_0$ for some $g_0$ which is holomorphic at $\infty$ and with $g_0(\infty) \neq 0$.\par
%\medskip
\begin{definition}
 Given  $A \in \cC^\epsilon(\partial \D, \GL(N, \bC))$, with $\epsilon \in ]0,1[$, consider  the problem 
consisting of  finding a continuous map $\Psi^+: \overline \D \to \bC^N$, holomorphic on $\D$, and a continuous map 
$\Psi^-: \bC \setminus \D \to \bC^N$, holomorphic on $\bC \setminus \overline \D$ and with at most a pole at $\infty$, so that 
\beq \label{problem} \Psi^+(\z) = A(\z) \cdot \Psi^-(\z)\ ,\qquad \z \in \partial \D\ .\eeq
A  {\it canonical system of $A$} is any collection of solutions $\Phi_j = (\Phi_j^+, \Phi_j^-)$, $1 \leq j \leq N$, 
of the problem (\ref{problem}) so that 
\begin{itemize}
\item[i)] $\Phi^+(\z) = [\Phi^+_1(\z), \dots, \Phi^+_N(\z)]$ is in $\GL(N, \bC)$ for any $\z \in \overline \D$; 
\item[ii)] $\Phi^-(\z) = [\Phi^-_1(\z), \dots, \Phi^-_N(\z)]$ is in $\GL(N, \bC)$ for any $\z \in \bC \setminus \D$; 
\item[iii)] the  order  $k$ of $\det \Phi^- $ at $\infty$ is equal to the sum of the  orders $k_j$ of the columns $\Phi^-_j$.   
\end{itemize}
If $\{\Phi_j = (\Phi_j^+, \Phi_j^-)\}$ is a canonical system of $A$, the orders $k_j$ of the  $\Phi_j^-$'s are called 
{\it partial indices of $A$\/}. The sum   $k = \sum k_j$ is called {\it total index of $A$\/}. \par
\end{definition}
An important fact is that, up to reordering,  {\it the partial indices and the total index depend only on $A$ and not 
on the considered   canonical system\/}.
We may now recall the following lemma  
by  Globevnik,  which can be considered as a corollary of N. P. Vekua's factorization theorem (\cite{Vek}).\par
\begin{lem} \label{factorization} (\cite{Gl}, Lemma 5.1)
Let $L\in \cC^\epsilon(\partial \D, \GL(N, \bC))$, with $\epsilon \in ]0,1[$. Then there is 
a map $\Theta: \bar \Delta \to \GL(N,\bC)$ in  $H^\epsilon(\overline \Delta, \bC^{N^2})$,  
such that 
\beq \label{Globevnik} 
L(\zeta)\cdot \overline{L(\zeta)^{-1}} = 
\Theta(\zeta) \cdot \Lambda(\zeta) \cdot \bar \Theta^{-1}(\zeta) 
\ \text{with}\  \Lambda(\zeta) = \left(\smallmatrix
\zeta^{k_1} & 0  & \dots & 0\\
0 & \zeta^{k_2}  & \dots & 0\\
\vdots & \vdots  & \ddots & \vdots \\
%0 & \dots  & \zeta^{k_{N-1}}& 0\\
0 & 0 & \dots  & \zeta^{k_N}
\endsmallmatrix\right)\eeq
 for any $\zeta \in \partial \Delta$, where  $k_1, \dots, k_N$ are  the partial indices of $A(\cdot) \=  \left. L(\cdot)\cdot \overline{L(\cdot)^{-1}}\right|_{\partial \D}$.
\end{lem}
%\bigskip
The integers $k_i$ of the previous lemma and the sum $k = \sum_{i=1}^N k_i$ are the same for all  maps  $L'  = M|_{\partial \D}\cdot L$ with $M: \bar \Delta \to \GL(N,\bC)$ in $H^\epsilon(\overline \Delta, \bC^{N^2})$. 
They are called {\it partial 
indices\/} and  {\it total\/}   {\it index\/} of $L$, respectively.\par
\smallskip
Consider now the map $G(\z)$ in (\ref{matrixG}) and  let  $\Theta^G$ be a map that gives a decomposition (\ref{Globevnik})  for $L(\z) = G^{-1}(\z)$. 
We set 
\beq \upG A \= \left({\upG \Theta}\right)^{-1} \cdot A \cdot \upG \Theta\ , \qquad \upG B \= \left({\upG \Theta}\right)^{-1} \cdot B \cdot \overline{\upG \Theta}\ .\eeq
It is immediate to realize  that
the linear map $\wh h \longmapsto \wt h =  \left({\upG \Theta}\right)^{-1} \cdot \wh h$ is 
  an isomorphism between  
 $\ker R_{A,B,G}$  the  space of 
solutions   of the problem
\beq\label{diagonalized} \left\{ 
\begin{array}{llll}
\bar \partial \wt h + \upG A\cdot \wt h + \upG B\cdot \overline{\wt h} = 0\ , &  \zeta \in \Delta \\
& & &\\
\wt h^i(\zeta) =  \zeta^{k_i}  \overline{\wt h^i}(\zeta)\ ,\ \ 1 \leq i \leq 2n, \ \ \   & \zeta \in \partial \Delta
\end{array}\right. \eeq
where the $k_i$ are the partial indices of  $L = G^{-1}$. \par
%\medskip
\begin{lem} \label{globevnik'slemma} The operator $R_{A,B,G}$  is surjective if and only if 
$  \dim \ker R_{A,B,G} = 2n  +k$, with $k =  \sum_{i = 1}^{2n} k_i$. Moreover, when  $A$ $= B$ $= 0$, $R_{0,0,G}$  is surjective
 if and only if $k_i \geq -1$ for any $1 \leq i \leq 2n$. 
 \end{lem}
\begin{pf} The first claim follows from Theorem \ref{indexR} and from 
\begin{align*}
 \dim \ker R_{A,B,G} & = 2n  - \frac{1}{i \pi} \int_{\partial \Delta} d\arg(\det(G)) \\
 &= 2n + 
\frac{1}{2 \pi i} \int_{\partial \Delta} d\arg\left(\det (G^{-1} \cdot \overline{G}) \right) \\
& =  2n +
\frac{1}{\pi i} \int_{\partial \Delta} d\arg(\det(\upG \Theta)) + 
\sum_{i = 1}^{2n}\frac{1}{2 \pi i} \int_{\partial \Delta} d\arg(\zeta^{k_i})\\ &=   2n +  k
\end{align*}
%
%
%$$ \dim \ker R_{A,B,G} = 2n  - \frac{1}{i \pi} \int_{\partial \Delta} d\arg(\det(G)) = $$
%$$ = 2n + 
%\frac{1}{2 \pi i} \int_{\partial \Delta} d\arg\left(\det (G^{-1} \cdot \overline{G}) \right) = $$
%$$ = 2n +
%\frac{1}{\pi i} \int_{\partial \Delta} d\arg(\det(\upG \Theta)) + 
%\sum_{i = 1}^{2n}\frac{1}{2 \pi i} \int_{\partial \Delta} d\arg(\zeta^{k_i})=   2n +  k\ , $$
where we used the fact  that  $\det(\upG \Theta)$ is holomorphic and never zero in $\D$. \par
%\smallskip
Assume now that  $A = B = 0$ and recall that the elements of $\ker R_{0,0,G}$ are in natural correspondence with 
the elements  $\wt h= (\wt h^1, \dots, \wt h^{2n}) \in H^\varepsilon(\D, \bC^{2n})$  that solve (\ref{diagonalized}) and hence  of the form
$\wt h^i(\zeta)  = \sum_{\ell \geq 0} a^i_\ell \z^\ell$ 
with coefficients   $a^i_\ell \in \bC$ so that the boundary conditions are satisfied, i.e. 
\beq \label{finedim}
\left\{\begin{array}{lll} a_\ell^i = 0 & & \hbox{when}\ \ell \geq \max\{k_i + 1, 0\}\\
\ \\
a_\ell^i =   \overline{a_{-\ell+ k_i}^i}   & & \hbox{when}\  k_i \geq 0 \ \ \text{and}\ 0 \leq \ell \leq k_i\ .
\end{array}\right.\eeq
From this, a simple check shows that
 $ \dim \ker R_{0,0,G} =$ $ \sum_{k_i\geq 0} (k_i + 1)$.  
 Since 
 $$ 2n +  k =  2n  - \sum_{k_i\leq -1}( |k_i | - 1 ) -( \#\{k_i\leq -1\} )  + \sum_{k_i\geq 0} k_i  = $$
$$=  (\# \{\ k_i \geq 0\}) + \sum_{k_i\geq 0} k_i  - \sum_{k_i\leq -1}( |k_i | - 1 ) =  \sum_{k_i\geq 0} (k_i + 1)  - \sum_{k_i\leq -1}( |k_i | - 1 )$$
 it follows that 
 $ \dim \ker R_{0,0,G} = 2n  + k$ if and only if  $ \sum_{k_i\leq -1}( |k_i | - 1 )  =  0$, i.e. $k_i \geq - 1$ for any $1 \leq i\leq 2n$. \end{pf}
\bigskip
 \subsection{The operator $\gR = (R_{0,0,G}, \gR_3, \gR_4, \gR_5)$ for convex domains in $\bC^n$}
 \begin{theo} \label{main2}ÊLet $D$  be a  domain in $(\bC^n, \Jst)$, 
with smooth boundary  and let $f_o :\overline  \Delta \to \overline D$ a stationary disk  $D$. If  there is  a neighborhood 
$\cU$ of $f_o(\overline D)$   where  $\cU \cap \overline D$ is  strictly linearly convex, then $\partial D$ is  good  for  $(\Jst, f_o)$.
\end{theo}
\begin{pf} We first  need the following:\par
\begin{lem}[\cite{Pa}, Prop. 2.36, Thm. 2.45]  Let $f_o:\overline  \Delta \to \overline D$ as  above. Then there exists a system 
of complex coordinates $(z^i)$ and a defining function $\rho$ for $\partial D$ on a neighborhood
$\cV$ of $f(\overline \D)$, such that  $f_o(\zeta) = (\zeta, 0, \dots, 0)$ and 
\beq\label{normalformbis}  \rho = - 1 + |z^1|^2 + \sum_{\a, \b = 2}^n \d_{\a \b} z^\a \overline{z^\b} + 
\Re\left(\sum_{\a, \b = 1}^n B_{\a \b} z^\a z^\b\right) + r(z^1, \dots, z^n)\eeq
with $r$ smooth function  so that $|r(z)| \leq c |z|^3$ for some $c >0$ 
for all $z  \in\cV$.
\end{lem}
Secondly, we need  the following  lemma, from which the theorem will follows almost immediately. There, we denote by $(z^i)$ the coordinates  in previous lemma  and  by
$ (z^i, w_i)$ the associated complex coordinates for $T^* \bC^n$ (see \S \ref{notation}). 
%\bigskip
\begin{lem} \label{indices} Let $\gR = (R_{0,0,G}, \gR_3, \gR_4, \gR_5)$ be the 
linear operator defined in  \S \ref{condcoer} using the coordinates $(z^i, w_j)$. Then: 
\begin{itemize}
\item[i)] The partial indices of  $G^{-1}$ are $k_1 = 2$, $k_2 = 0$ and $k_j = 1$ for all $j \geq 2$. In particular,
$R_{0,0,G}$ is surjective and $\dim \ker R_{0,0,G} = 4n +1$.
\item[ii)] The restrictions of  $\gR_3$, $\gR_4$  on $\ker R_{0,0,G}$ are surjective.
\end{itemize}
\end{lem}
\begin{pf} (i) If $\rho$ is the defining function (\ref{normalformbis}),  the components of the  function 
$\varrho(\zeta, \a)$,  defined in (\ref{varrhozeta}),  are  (up to multiplication by a nowhere vanishing smooth function)
$$\varrho^1 =  - 1 + |z^1|^2 + \sum_{\a, \b = 2}^n \d_{\a \b} z^\a \overline{z^\b} + 
\Re\left(\sum_{\a, \b = 1}^n B_{\a \b} z^\a z^\b\right) + O(|z|^3)$$
$$\varrho^2 = i \left\{2 |z^1|^2 (\z^{-1} w_1 - \overline{\z^{-1}}\overline{w_1}) - (z^1 \zeta^{-1} w_1 + \overline{z^1} \overline{\z^{-1}} \overline{w_1} )(\overline{z^1} - z^1)\right\} + O(|z|^2)$$
$$\varrho^{2 \a - 1} = 2 |z^1|^2 (\z^{-1}w_\a + \overline{\z^{-1}} \overline{w_\a}) - \phantom{aaaaaaaaaaaaaaaaaaaaaaaaaaaaaaaa}$$
$$ -  (z^1 \zeta^{-1} w_1 + \overline{z^1} \overline{\z^{-1}} \overline{w_1} )\left\{(\d_{\a \b} \overline{z^\b} + B_{\a \b} z^b) + (\d_{\a \b} z^\b + \overline{B_{\a \b}} \overline{z^\b}) \right\} + O(|z|^2)$$
$$\varrho^{2 \a} = i \left\{2 |z^1|^2 (\z^{-1} w_\a - \overline{\z^{-1}} \overline{w_\a}) - \right. \phantom{aaaaaaaaaaaaaaaaaaaaaaaaaaaaaa}$$
$$ \left. -  (z^1 \zeta^{-1} w_1 + \overline{z^1} \overline{\z^{-1}} \overline{w_1} ) \left\{(\d_{\a \b} \overline{z^\b} + B_{\a \b} z^b) - (\d_{\a \b} z^\b + \overline{B_{\a \b}} \overline{z^\b}) \right\}\right\} + O(|z|^2)$$
with $2 \leq \a \leq n$. Hence, the matrix  (\ref{matrixG}) is (up to reordering of columns) 
$$G(\zeta) = \left(\begin{matrix}ÿG_1(\zeta) & 0\\ 0 & G_2(\zeta)\end{matrix} \right)\qquad \text{with}$$
$$G_1(\zeta) = \left.\left(\begin{matrix} \frac{\partial \varrho^1}{\partial z^1} &  \frac{\partial \varrho^1}{\partial w_1}\\
\frac{\partial \varrho^2}{\partial z^1} &  \frac{\partial \varrho^2}{\partial w_1}\end{matrix}\right)\right|_{\wh f_o(\partial \D)} = 
\left(\begin{matrix} \z^{-1} & 0\\
i (\z^2 + 1) \z^{-2} & i (\z^2 + 1) \z^{-1}\end{matrix}\right)$$
$$G_2(\zeta) = \left.\left(\begin{matrix} \frac{\partial \varrho^{2\a -1}}{\partial z^\b} &   \frac{\partial  \varrho^{2\a -1}}{\partial w_\g}\\
\frac{\partial \varrho^{2\a}}{\partial z^\b} &   \frac{\partial \varrho^{2\a}}{\partial w_\g}\end{matrix}\right)\right|_{\wh f_o(\partial \D)} = 
\left(\begin{matrix}  
 -  2  B_{\a \b} + \d_{\a \b}  &  2  \z^{-1} \\
 -  2 i B_{\a \b} z^b + i \d_{\a \b} z^\b &  i 2 \z^{-1} \end{matrix}\right)$$
Under the assumption that the real Hessian $H(\rho)_{ij}$  is positive definite at all points of  $f_o(\partial \D)$, the partial indices of the matrix  $G^{-1}_2(\zeta)$ are known  to be all equal to $1$. A complete proof of this can be found in \cite{ST},  Lemma 3.2,  being   $G_2(\zeta)$  equal to the lower right block  of the matrix in (3.10) of \cite{ST}. \par
For what concerns the block $G^{-1}_1(\zeta)$, notice that  for any $\z \in \partial \D$ one has that $A(\z) = G^{-1}_1(\z) \cdot \overline{G_1(\z)}  = \left(\begin{matrix}\z^2 & 0 \\ - 2 \z & -1\end{matrix}\right)$ and hence $A$ admits  the columns of 
$$\Phi^{+}(\zeta) = \left(\begin{matrix}1 & 0 \\ 0 & 1\end{matrix}\right)\ ,\qquad \Phi^{-}(\zeta) = \left(\begin{matrix} \frac{1}{\z^2} & 0 \\  \frac{2}{\z} & -1\end{matrix}\right)$$
as canonical system. 
Hence, by  Lemma \ref{globevnik'slemma}, we conclude that  $k_1 = 2$ and $k_2 = 0$, since these  are  the orders of the columns of $\Phi^{-}$.\par
\smallskip
\noindent (ii) Let $I_{x_o} \subset T_{x_o} D$ be the indicatrix of the Kobayashi metric of $D$ at $x_o$. By Thm. 4.8 in \cite{Pa}, there exists a neighborhood $\cW \subset D$ of $x_o$ and a neighborhood $\cW' \subset \partial I_{x_o}$ of $v_o$, so that 
for any $x \in \cW$,  $v \in \cW' $
 there exists exactly two stationary disks $f^{(x, v_o)}, f^{(x_o, v)}: \overline{\D} \to \overline{D}$ 
satisfying
$$f^{(x, v_o)}(0) = x\ ,\ f^{(x, v_o)}_*\left(\left.\frac{\partial}{\partial x}\right|_0\right) = v_o\ ,\  f^{(x_o, v)}(0) = x_o\ ,\ \ f^{(x_o, v)}_*\left(\left.\frac{\partial}{\partial x}\right|_0\right) = v$$
For each of them,   there is a  unique stationary lift $\wh f^{(x, v_o)}$ and $\wh f^{(x_o, v)}$ satisfying  certain normalizing conditions  (i.e. so that   $\z \cdot \wh f^{(x, v_o)}(\z)$ and 
 $\z \cdot \wh f^{(x, v_o)}(\z)$ are  the so called {\it dual maps\/} - see \cite{Pa}, Def. 2.10). These lifts 
depends smoothly on the coordinates of the point $x$ and the vector $v$ and for any  curves $\g_t \in D$ and $\g'_t \in T_{x_o} D$ with
$\g_0 = x_o$ and $\g'_0 = v_o$, the 1-parameter families of stationary lifts $\wh f_t \= \wh f^{(\gamma_t, v_o)}$ and $\wh f'_t \= \wh f^{(x_o, \gamma'_t)}$
are so that $\wh h(\zeta) = \left.\frac{d \wh f_t(\z)}{dt}\right|_{t= 0}$ and $\wh h'(\zeta) = \left.\frac{d \wh f_t'(\z)}{dt}\right|_{t= 0}$ are  in $ \ker R_{0,0,G}$. Moreover, by construction, 
$$\gR_3(\wh h) = \pi(\wh h(0)) = \dot \g_0 \in \bC^n\ ,\ \ \gR_4(\wh h, \s) = \left.\frac{\partial (\pi \circ \wh h)}{\partial x}\right|_0  - \s v_o =  \dot \g'_0  - \s v_o \in \bC^n\ .$$
Since $v_o$ is transversal to $\partial I_{x_o}$, by the arbitrariness of  $\g_t$ and $\g' \in \partial I_{x_o}$ it follows that 
$\gR_3|_{\ker R_{0,0,G}}$ and $\gR_4|_{\ker R_{0,0,G}}$ are both surjective.
\end{pf}
By the previous lemma, $\dim \ker R_{0,0,G} \cap \ker \gR_3 \cap \ker \gR_4 = 1$. So, in order to conclude,   we only need to check that $\gR_5|_{\ker R_{0,0,G} \cap \ker \gR_3 \cap \ker \gR_4}$ is surjective onto $\bR$ or, equivalently, that there is  $0 \neq \wh h \in   \ker R_{0,0,G} \cap \ker \gR_3 \cap \ker \gR_4$ so that  
$\gR_5(\wh h) = \wh h_1(1, 0, \dots, 0) \neq 0$. But  an element of this kind is given by $\wh h(\zeta) = \left.\frac{d (\varphi_t(\wh f_o(\z)))}{dt}\right|_0 = (\zeta, 1, 0, \dots 0)$,  
where  we denote by $\varphi_t$  the diffeomorphism considered in Lemma \ref{littlelemma}, and the proof is concluded.
\end{pf}
\begin{rem} Lemma \ref{indices} (i) corrects and generalizes a computation in \cite{CGS}, where, by a minor mistake,  the   partial indices of $G^{-1}$ in case  $D = B^n$  are claimed to be  all equal to $1$. 
\end{rem}
\bigskip
\section{Other non-singular foliations by stationary disks}
\subsection{Foliations of horospherical type}
As before, $(M, J)$ is an almost complex manifold of dimension $2n$. 
Let   $x_o \in M$  and consider a  Riemannian metric $<, >$  on a neighborhood $\cU$ so that $<,>|_{x_o}$ is $J$-Hermitian. For instance, if  $\cU$ is identified with an open subset of $\bC^n$ so that 
$J|_{x_o} = \Jst|_{x_o}$,   we may assume that  $<,>$  is  the standard Hermitian metric of $\bC^n$. Denote also by  $\n$  the Levi-Civita connection of $<,>$. \par
\begin{definition} Let   $f: \overline \D \to M $ be a $J$-holomorphic  disk, which is  $\cC^1$ up to the boundary and with  $v _o = f_*\left( \left.\frac{\partial}{\partial x}\right|_1\right) \neq 0$.  We call   {\it parameter  of tangency at $x_o = f(1)$\/}  the real number
\beq p(f; x_o) \=  \left< \n_{v_o} \left(f_*\left(\frac{\partial}{\partial x}\right)\right), J v_o\right> \ . \eeq
\end{definition}
%\medskip
This number  depends on the first order jet of   $<,>$ at $x_o$, but if two $J$-holomorphic disks $f$, $h$ 
are so that  
$$x_o = f(1) = h(1)\ , \qquad v _o = f_*\left( \left.\frac{\partial}{\partial x}\right|_1\right) = h_*\left( \left.\frac{\partial}{\partial x}\right|_1\right)\ ,$$
then their  parameters of tangency are the same for a choice of  $<,>$ if and only if they are the same for any other choice of  the  metric.  In fact, if we consider a new metric  $<,>'$ with Levi-Civita connection $\n'$,  then $S = \n'-\n$ is a tensor field of type $(1,2)$ so that
 $$\left.\left(\n_{v_o} \left(f_*\left(\frac{\partial}{\partial x}\right)\right) -  \n_{v_o} \left(h_*\left(\frac{\partial}{\partial x}\right)\right)\right)\right|_1 = S(v_o, v_o - v_o) = 0$$
Moreover,   a simple computation shows that  {\it any  disk  $h = f \circ \varphi$ where $\varphi \in \Aut(\D)$ with $\varphi(1) = 1$, $\varphi'(1) = 1$,  satisfies\/}
$$\n_{v_o} \left(h_*\left(\frac{\partial}{\partial x}\right)\right) = \lambda J f_*\left(\left.\frac{\partial}{\partial x}\right|_1\right) = \lambda J v_o \qquad \text{for some}\ \  \l \in \bR\ .$$
Therefore,  for any given $\wt \l$, one can choose $\varphi$ so that  $p(f \circ \varphi; x_o) = \wt \l$. Moreover,   $p(f \circ \varphi; x_o) = p(f; x_o)$ if and only if  $\varphi = Id_\D$ and $f = h$. \par
\smallskip
Consider now  a bounded, strictly convex domain in  $(\bC^n, \Jst)$ with smooth boundary and let $x_o \in \partial D$ and $\nu$  the outward unit normal to $\partial D$ in $x_o$. 
By  \cite{CHL}, Thm. 2,   for any $v_o \in T_{x_o}M$ so that $<\nu, v_o> > 0$  and for any $\lambda \in \bR$,  there exists a unique stationary disk $f^{(v_o, \lambda)}: \overline \D \to \overline D$ so that
\beq f^{(v_o, \lambda)}(1) = x_o\ , \ \  f_*\left( \left.\frac{\partial}{\partial x}\right|_1\right) = v_o\ ,\ \ \ p(f; x_o) = \l\ .\eeq
If we denote by $H_{x_o} = \{ \ v \in T_{x_o} \bC^n\ : \ < \nu, v> > 0\ \} \subset T_{x_o} \bC^n$, we have that 
 the exponential map 
 $\Phi^{(D, x_o)}: H_{x_o} \times (\D \setminus \{1\})\to \overline{D} \setminus \{x_o\}$ defined by  $\Phi^{(D, x_o)}(v; \z) = f^{(v, 0)}(\z)$
 is a diffeomorphism. \par
 We now consider  the following definition. As before, $D$ is a smoothly bounded, strictly pseudoconvex domain in the almost complex manifold $(M,J)$ and,  for any given $x_o \in \partial D$,  we denote by $\nu$  the outward unit normal to $\partial D$ in $x_o$ w.r.t. some Riemannian metric $<,>$, which is $J$ Hermitian at $x_o$. Finally,  for any real number   $a > 0$, we denote by $\cC^{(a)}$ the open cone 
$$\cC^{(a)} = \{\ v \in T_{x_o}ÊM\ : \ < v, \nu> > a\ \} \subset T_{x_o} M\ . $$
\begin{definition}\label{horo} For any $x_o \in \partial D$ and $a >0$,  let   $\cG^{(x_o)}$ be the family  of stationary disks $f: \overline \D \to \overline D$ with  $f(1) = x_o$  and by  $\cG^{(x_o, a)} \subset \cG^{(x_o)}$ the subfamily of disks with $ f_*\left( \left.\frac{\partial}{\partial x}\right|_1\right)  \in \cC^{(a)}$. Denote also by  $D^{(x_o, a)} \subset D$  the union of all images of the disks in $\cG^{(x_o, a)}$.\par
We say that  $\cG^{(x_o)}$  {\it  is a foliation  of  horospherical type for $D$\/} (resp. {\it  $\cG^{(x_o, a)}$ is a good foliation for $D^{(x_o, a)}$\/}) if  the following conditions are satisfied:
\begin{itemize}
\item[i)] for any $v \in T_{x_o} M$ so that $<v, \nu> >0$ (resp. for any $v \in \cC^{(a)}$) and for any $\l \in \bR$  there exists 
a unique $f^{(v, \l)} \in \cG^{(x_o)}$ so that
\beq  f^{(v, \l)}_*\left( \left.\frac{\partial}{\partial x}\right|_1\right) = v\ ,\ \ \ p(f; x_o) = \l\eeq 
\item[ii)]   the map   $\exp: \overline {B^n} \setminus \{y_o\} \to \overline D \setminus\{x_o\} $, $y_o \= (1, 0, \dots, 0)$,   defined by 
\beq \exp(\Phi^{(B^n, y_o)}(v, \zeta)) \= f^{(v, 0)}(\zeta)\eeq
is a diffeomorphism on $B^n$ (resp. on $B^n{}^{(y_o, a)}$), extends  smoothly at all points of the closure, different from $y_o$,  and induces an homeomorphism between the closures of  the two domains.
%
%previous an identification  $(T_{x_o} M, J_{x_o}) \simeq (\bC^n, \Jst)$,  
% the map
%\beq  \exp: \wt B^n  \subset \wt \bC^n \longrightarrow  \wt D \ \ , \qquad  \ \exp(v, [v]) \= \wt{f^{(v)}}(|v|)\ ,\eeq
%between the blow up $\wt B^n$ at $0$ of   $B^n\subset \bC^n $  and $\wt D$ 
%is smooth and with a smooth extension  up to the boundary, which induces a diffeomorphism  between the boundaries  $\exp|_{\partial B^n}: \partial B^n \to \partial D$.
 \end{itemize}
If  $\cG^{(x_o, a)}$ with $a> 0$ is  a good foliation  for $D^{(x_o, a)}$ ,  we say that     {\it  $D^{(x_o, a)} \subset D$ is a good conical subdomain  with vertex in $x_o$}. If  $\cG^{(x_o)}$  is  a foliation of horospherical type, we say that $x_o$ is a {\it center at infinity for $D$\/} and  {\it $D$ is of horospherical type\/}.\end{definition}
By the results in \cite{CHL},  any strictly convex domain  $D$ in $(\bC^n, \Jst)$ is a domain of horospherical type with  center at infinity at any point of the boundary. \par
%\smallskip
%\bigskip
\subsubsection{Stability of foliations of horospherical type }
In analogy with  \S \ref{stability},  let us consider the nonlinear operator $\cR'_{(J_o, x_o, v_o, \nu)} = (\cR'_1, \dots, \cR'_6)$ from 
$\cC^{\alpha, \varepsilon}(\overline\Delta; \bC^{2n})  \times  \cC^{\epsilon}(\partial \Delta; \bR) $
into $ 
\cC^{\alpha-1, \varepsilon}(\overline\Delta; \bC^{2n})  \times   \cC^{\epsilon}(\partial \Delta; \bR^{2n+1}) \times \partial D\times \bC^n \times \bR \times \bR$, where
% $\cR'_1$, $\cR'_2$ and $\cR'_6$  coincide with the components $\cR_1$, $\cR_2$ and $\cR_5$ of  (\ref{nonlinearoperator}), while $\cR_j$, $j = 3, 4,5$, are defined by
$$\cR'_1(\wh f, \l) = \overline\partial_{\bJ_o}\wh f\ \ ,\ \ 
\cR'_2(\wh f, \l ) = r(\zeta,  \lambda(\zeta), \wh f(\zeta))\ ,$$
$$ \cR'_3(\wh f, \l) = \pi(\wh f)|_{\zeta = 1} - x_o\ , \  \cR'_4(\wh f, \l) = \pi(\wh f)_*\left(\left.\frac{\partial}{\partial x} \right|_{\zeta = 1}\right) -  v_o\ ,$$
$$  \cR'_5(\wh f, \l)=  p(\pi(\wh f); x_o) - \nu \ ,$$
$$  \cR'_6(\wh f, \l) =  \wh f \left(\pi(\wh f)_*\left(\left.\frac{\partial}{\partial x} \right|_{1}\right)\right) - 1\ .$$
Given a stationary disk $f_o: \overline\D \to \overline D$   with $x_o = f(1)$, $v_o = f_o{}_*\left(\left.\frac{\partial}{\partial x} \right|_{\zeta = 0}\right)$ and $p(f_o; x_o) = \nu_o$,  we say that    $\partial D$ is {\it  a horospherically good boundary for  $(J_o, f_o)$\/}  if $f_o$ admits  a lift $\wh f_o$  so that    $(\wh f_o, \l_o)$ is a solution of  $ \cR'_{(J_o, x_o, v_o, \nu_o)}(\wh f, \l) = 0$   and the  tangent  operator $\gR'$ of $ \cR'_{(J_o, x_o, v_o, \nu_o)}$ at $(\wh f_o, \l)$ is invertible.\par
Again, by the Implicit Function Theorem,   if $\partial D$ is a horospherically good boundary  for   $(J_o, f_o)$,   there is  a neighborhood $\cV \subset \partial D$ of $x_o$,  a neighborhood $\cW \subset T D$ of $v_o$, with $\pi(\cW) = \cV $ and a real number  $\varepsilon > 0$ so that,  for any $x \in \cV$, $v\in \cW$, $|\nu - \nu_o| < \varepsilon$    and    $\| J - J_o\|^{(1)}_{\overline D} < \varepsilon$,  there exists a unique disk $f$ in $D$ with 
\beq  f(1) = x\ ,\qquad f_*\left(\left.\frac{\partial}{\partial x}\right|_{\zeta = 1}\right) = v\ \ \ \text{and }\ \ p(f; x) = \nu\ ,\eeq
which is stationary  for $D$ w.r.t. the almost complex structure $J$. 
The dependence of  $f$ on  $x$,   $v$, $\nu$  and   $J$ is differentiable and,  given  $m_o > 0$ and a metric $g$, one can choose $\varepsilon$, $\cW$ and $\cV = \pi(\cW)$, so that  $\sup_{\zeta \in \overline \D} \dist_g(f(\zeta), f_o(\zeta))  < m_o$. So, in analogy with Proposition \ref{main3}, we have: 
\begin{prop} \label{giorgio}
Let  $D^{(x_o, a)} \subset D$, $a > 0$,  be a good conical subdomain    w.r.t. to $J_o$ with vertex in $x_o  \in \partial D$. If   $\partial D$ is a good boundary for $(J_o, f_o)$ for  any stationary disk  $f_o \in \cG^{(x_o, a)}$, there exists  $\varepsilon > 0$ and an open neighborhood $\cU \subset \partial D$ of $x_o$  so that  for any $J$ with   $\| J - J_o\|^{(1)}_{\overline D} < \varepsilon$ and any $x \in \cU$ and $|a'- a| < \varepsilon$, the point $x$  is   vertex  for a good  foliation  for  $D^{(x, a')}$ relatively to  the almost complex structure $J$.
\end{prop}
\begin{pf} The proof  can be obtained following the same 
steps of  the proof of Prop. 6  in \cite{CGS} and we give here only a sketch of it.  First of all, using the Implicit Function Theorem and the compactness of $\overline{ \cC^{(a)} }\cap S^{2n-1} \subset T_{x_o} M$,  one can determine $\cU$ and $\varepsilon$ so that $\cG^{(x, a')}$  satisfies (i) for Definition \ref{horo} for any almost complex structure such that  $\| J - J_o\|^{(1)}_{\overline D} < \varepsilon$ and for any  $x \in \cU$,  $|a'- a| < \varepsilon$. Using the Implicit Function Theorem once again, one can also assume that for all these $J$, $x$ and $a$, the map ``$\exp$'',  defined in (ii) of that definition,  is a local diffeomorphism at all points. It remains to be checked that  
$\cU$ and $\varepsilon$ can be chosen so that ``$\exp$'' is also injective.  From this and  a possible further restriction of $\cU$ and $\varepsilon$, we  obtain that ``$\exp$'' is a diffeomorphism and satisfies all other requirements of (ii). To prove injectivity,  one may argue by contradiction as in Step 2 of the proof of Prop. 6 in \cite{CGS}. In fact, if one assumes that ``$\exp$'' is never injective for any choice of $\cU$ and $\varepsilon$, one can construct sequences of complex structures $J_j$, of vertices $x_j$ and of pairs  $y_j \neq y'_j \in B^n$, so that $J_j \to J_o$, $x_j \to x_o$ and corresponding
exponential maps $\exp^{(j)}$ are so that $\exp^{(j)}(y_j) =\exp^{(j)}(y'_j)$ for all $j$. Using compactness and Implicit Function Theorem, one can select a subsequence
$z_{j_m} = \exp^{(j)}(y_{j_m})$, with   $y_{j_m} \to y_o$, $y'_{j_m} \to y'_o$ with $y_o \neq y'_o$ and $z_{j_m}  \to z_o = \exp(y_o) = \exp(y'_o) \in \overline{D^{(x_o, a)}}$, contradicting the hypothesis of bijectivity of  ``$\exp$'' on $\overline{D^{(x_o, a)}}$. 
\end{pf}
We remark that the tangent operator $\gR' = (\gR'_1,  \dots, \gR'_6)$ of $ \cR'_{(J_o, x_o, v_o, \nu_o)}$ at $(\wh f_o, \l)$ is so that 
$(\gR'_1, \gR'_2) = R_{A, B, G}$ (see \S \ref{R_{A,B,G}} for definition)  and hence it  coincides with  operator $R_{0,0,G}$ when  $D \subset \bC^n$.  If $D$ is a strictly linearly convex 
domain  in $\bC^n$,  by Lemma \ref{indices} (1), the dimension of $\ker (\gR'_1, \gR'_2) = 4n + 1$. From this,  the results in \cite{CHL} and  a line of argument which is essentially the same of  the proofs of Lemma \ref{indices} (1) and Theorem \ref{main2}, one gets that  $\gR' = (\gR'_1, \dots, \gR'_6)$ is invertible also  in this case.  By Proposition \ref{giorgio},   the following result is obtained.\par
\begin{theo} Let $D \subset M$ be a  smoothly bounded, strongly  pseudoconvex 
domain in an almost complex manifold $(M, J_o)$ and $a > 0$ be any fixed  positive real number.  If there is a local diffeomorphism $\varphi: \cU \subset M \to \bC^n$, so that  $\wh D = \varphi(D)$ is  a strictly linearly convex domain $\wh D \subset \bC^n$   and $\varphi_*(J_o)$ is sufficiently close to $\Jst$ in a  $\cC^1$-norm, then, for any $x_o \in \partial D$, the subset  $D^{(x_o, a)} \subset D$   is a good conical subdomain. 
\end{theo}
%\bigskip

\end{document}